\documentclass[a4paper,10pt]{article}

\usepackage[latin1]{inputenc}
\usepackage{bbm} 
\usepackage{mathrsfs} 

\usepackage{amsfonts}
\usepackage{amsmath}
\usepackage{amssymb}
\usepackage{amsthm}
\usepackage[american]{babel}
\usepackage[T1]{fontenc}
\usepackage[dvips]{graphicx}

\newtheorem{Theorem}{Theorem}[section]
\newtheorem{Definition}[Theorem]{Definition}

\newtheorem{Lemma}[Theorem]{Lemma}
\newtheorem{Corollary}[Theorem]{Corollary}

\newtheorem{Remark}[Theorem]{Remark}
\newtheorem{ass}[Theorem]{Assumptions}

\newcommand{\taunew}{T}

\newcommand{\ritardo}{R}

\newcommand{\ud}{\mathrm{d}}
\newcommand{\be}{\begin{equation}}
\newcommand{\ee}{\end{equation}}
\newcommand{\nd}{\stackrel{def}{=}}
\newcommand{\s}{\hspace{.05in}}

\newcommand{\lla}{\left\langle}
\newcommand{\rra}{\right\rangle}

\def\clip{C_{Lip}^1}
\def\vertv{\vert_{_{\V}}}
\def\t{\taunew}

\def\A{A_0}
\def\V{V^\prime}
\def\ranglev{\rangle_{V^\prime}}
\def\lipd{ ]\kern-1pt_{_{L}}}
\def\lips{[}
\def\H{\mathcal{H}}
\def\HAM{F}
\def\e{\varepsilon}

\begin{document}
\title{On the Dynamic Programming approach to economic models governed by DDE's}
\author{Giorgio Fabbri\footnote {Dipartimento di matematica \emph{Guido Castelnuovo},
universit\`a \emph{La Sapienza} Roma, e-mail:
fabbri@mat.uniroma1.it}, Silvia Faggian\footnote {Universit\`a
\emph{LUM - Jean Monnet,} Casamassima, Bari, e-mail:
faggian@lum.it}, Fausto Gozzi\footnote {Dipartimento di Scienze
Economiche ed Aziendali, Universit\`a \emph{LUISS - Guido Carli}
Roma, e-mail: fgozzi@luiss.it}}
\date{}
\maketitle

\begin{abstract}
In this paper we consider a family of optimal control problems for
economic models whose state variables are driven by Delay
Differential Equations (DDE's). We consider two main examples: an AK
model with vintage capital and an advertising model with delay
effect. These problems are very difficult to treat for three main
reasons: the presence of the DDE's, that makes them infinite
dimensional; the presence of state constraints; the presence of
delay in the control. Our main goal is to develop, at a first stage,
the Dynamic Programming approach for this family of problems. The
Dynamic Programming approach has been already used for similar
problems in cases when it is possible to write explicitly the value
function $V$ (see \cite{FGAK}). Here we deal with cases when the
explicit form of $V$ cannot be found, as most often occurs. We
carefully describe the basic setting and give some first results on
the solution of the Hamilton-Jacobi-Bellman (HJB) equation as a
first step to find optimal strategies in closed loop form.
\end{abstract}

\section{Introduction}





In this paper we want to develop the Dynamic Programming approach
for a family of optimal control problems related to economic models
governed by Delay Differential Equations (DDE's).

The presence of DDE's makes the problem difficult to treat. One
possible way of dealing with DDE's - the one we choose - is
rewriting the problem as an optimal control problem governed by
ODE's in a suitable Hilbert space. Although such infinite
dimensional optimal control problems have already been studied, the
present literature does not cover our case, as it does not include
the following features:
\begin{itemize}
    \item the presence of unbounded operators coming from the DDE which is
not analytic and does not satisfy smoothing assumptions;
    \item the presence of state/control constraints (which is indeed peculiar of
    economic models);
    \item the fact that the delay appears in the state and in the control (causing the control
operator to be possibly unbounded).
\end{itemize}
We stress the fact that these difficulties are the rule in economic
models governed by DDE's.

Here we consider problems with linear DDE's and \emph{concave}
objective functional: concavity will play a key role in the paper.
When concavity lacks, one can still apply Dynamic Programming in
the framework of viscosity solutions - which we avoid here.
Nevertheless, we address the reader to \cite{CILuserguide} for a
standard reference on viscosity solutions.

We remark that this is a first step in treating such kind of
problems. We already studied thoroughly in \cite{FGAK} a case where
explicit solution of the associated Hamilton-Jacobi-Bellman (HJB)
equation can be found (in such case the problem is much easier to
treat). Here we want to develop the Dynamic Programming approach in
those cases when explicit solutions of the associated HJB equation
are not available. We here develop the finite horizon case. The
infinite horizon case can be treated with our method using a
limiting procedure when the horizon goes to $+\infty$ but we leave
it for future work\footnote{In thus respect we can say that the
finite horizon case is as a first step towards the infinite horizon
one.}.

The main result of the paper is that the value function of the
problem is a solution, in a suitable weak sense, of the HJB
equation. This a first step towards the so-called Verification
Theorem which is a powerful tool to study the optimal paths of the
problem and which is the subject of our current research.

We concentrate on two main examples: an AK model with vintage
capital, taken from \cite{Boucekkine} (see also \cite{BoucekkineMPS}
and \cite{FGAK}) and an advertising model with delay effects (see
\cite{GMTrento,GMS}) that are exposed in Section
\ref{sezioneesempi}.

The plan of the paper is the following. In Section
\ref{sezioneesempi} we present the applied examples. In Section
\ref{sezioneDP} we recall the basic steps of the Dynamic Programming
approach and we give an overview of the current literature on the
Dynamic Programming for infinite dimensional optimal control
problems. In Section \ref{sezioneriscrittura} we rewrite the state
equation of such problems as an ODE in a suitable Hilbert space,
concentrating on the first example, as the second can be rephrased
similarly. In Section \ref{sectionHJB} we write the resulting
infinite dimensional optimal control problem and its HJB equation.
Section \ref{sezioneultraweak} we show our main result: the
existence of an ultraweak solution of the HJB equation. The Appendix
\ref{sezioneappendice} contains some definition and proof that may
be useful for the reader.

\section{Two examples}
\label{sezioneesempi} We present the two applied problems motivating
this paper.

\subsection{An AK model with vintage capital}

We consider here an optimal control problem related to a
generalization of the model presented by Boucekkine, Puch, Licandro
and Del Rio in \cite{Boucekkine}. Indeed, we assume that the system
is ruled by the same evolution law as the one
 in \cite{Boucekkine}, that is, an AK growth model with
a stratification on the capital. Besides, we consider the finite
horizon problem with a (more) general concave target functional,
as specified later. The analysis of such a  model proves
interesting in the study of short run fluctuations and of
transitional dynamics: the reader is referred to \cite{Boucekkine}
for a deep discussion upon this and other related matters. The
model of \cite{Boucekkine} is an infinite horizon model, while
here we consider the finite horizon case. As mentioned in the
introduction, this is a first step towards the infinite horizon
case.

The AK-growth model with vintage capital is based on the following
accumulation law for capital goods
\[
k(s)=\int_{s-\ritardo}^s i(\sigma) \ud \sigma
\]
where $i(\sigma)$ is the investment at time $\sigma$. That is,
capital goods are accumulated for the length of time $\ritardo $
(scrapping time) and then dismissed. Note that such an approach
introduces a differentiation in investments that depends on their
age.   If we assume a linear production function, that is
\[
y(s)=ak(s)
\]
where $y(s)$ is the output at time $s$ (note that "$AK$" reminds
of the linear dependence of the dynamic from the trajectory - a
constant $A$ multiplied by
 $K$; such constant $A$ is $a$ in our case), and we assume also
the accounting relation
\[
y(s)=c(s)+i(s),
\]
meaning that at every time the social planner chooses how to split
the production into consumption $c(s)$ and investment $i(s)$, then
 the state equation may be written into infinitesimal terms
as follows
\[
\dot{k}(s)= ak(s)- ak(s-\ritardo ) - c(s) + c(s-\ritardo ),\ \
s\in[t,\taunew]
\]
{\it i.e.} as a DDE. The time variable $s$ varies in $[t,\taunew]$,
with $t$ the \textit{initial time} and $\taunew$ the (finite)
\textit{horizon} of the problem.
Indeed, the social planner has to
maximize the following functional
\begin{equation}\label{target}
\int_{t}^\taunew e^{-\rho s}h_0(c(s)) \ud s + \phi_0(k(\taunew))
\end{equation}
where $h_0$ and $\phi_0$ are concave u.s.c. utility functions. We
recall that in \cite{Boucekkine} the horizon is infinite and
$\phi_0=0$. Moreover the instantaneous utility is CRRA (i.e. Costant
Relative Risk Aversion), that is the function $h_0$ is of type
$h_0(c)= \frac{c^{1-\sigma}}{1-\sigma}$, which satisfies our
assumptions as a subcase.

Observe that we take the starting time $t$ to be variable to apply
the finite horizon dynamic programming.

We assume that the capital at time $s$ (and consequently the
production) and the consumption at time $s$ cannot be negative:
\begin{equation}\label{eq:vincoliak}
    k(s)\ge 0, \qquad c(s) \ge 0,\quad  \forall s \in [t,\taunew]
\end{equation}
These constraints are different from the more restrictive and more
natural ones of \cite{Boucekkine}, where also the investment path
$i(\cdot) $ was assumed positive.

The main reason for such a choice is technical: we cannot apply the
strong solution approach that we use in this work with mixed
constraints such as those in \cite{Boucekkine}. The treatment of
mixed constraints is also left for future work. We mention indeed
that the optimal solutions for the problem without mixed constraints
may satisfy in some cases the positivity of investments, yielding
the solution also for the problem with mixed constraints.

In order to take the constraints into account, we assume that the
consumption (that is, the control variable of the system) lies in
the following admissible set
\[
{\mathcal{A}} \nd \{c(\cdot)\in L^2([t, \taunew],\mathbb{R})\; : \;
c(\cdot) \geq 0\; and \; k(\cdot)\geq 0\ a.e.\ in\  [t, \taunew]\}.
\]

\subsection{An advertising model with delay effects}

Another example of optimal control problems driven by DDE's is the
following a dynamic advertising model presented in the stochastic
case in the papers \cite{GMS,GMTrento}, and, in deterministic one,
in \cite{GozziFaggian} (see also \cite{hartladv} and the references
therein for related models)\footnote{We observe that also other
models of delay type arising in economic theory can be treated with
our tools (see \textit{e.g.} the paper by \cite{BoucekkineMPS}).}.

Let $t\ge 0$ be an initial time, and $T>t$ a terminal time
($T<+\infty $ here). Moreover let $\gamma(s)$, with $0\leq t\le s
\leq T$, represent the stock of advertising goodwill of the product
to be launched. Then the general model for the dynamics is given by
the following controlled Delay Differential Equation (DDE) with
delay $R>0$ where $z$ models the intensity of advertising spending:
\begin{equation}
\left\lbrace
\begin{array}{ll}
\dot \gamma(s)= a_{0} \gamma(s) +\int_{-R}^{0}
\gamma(s+\xi ) da_{1}(\xi ) +b_{0}z(s)+\int_{-R}^{0} z(s+\xi ) db_{1}(\xi ) \;\;\; s \in \lbrack t,\taunew] \\
\\
\gamma(t)=x;\quad \gamma(\xi )=\theta (\xi ),\;z(\xi )=\delta (\xi
)\;\;\forall \xi \in \lbrack t-R,t],
\end{array}
\right.
\label{eq:SDDE}
\end{equation}
with the following assumptions:

\begin{itemize}
\item   $a_{0}$ is a constant factor of image deterioration in absence of
advertising, $a_0 \leq 0$;

\item  $ a_{1}(\cdot )$ is the distribution of the forgetting
time, $a_1(\cdot) \in L^2([-R,0];\mathbb{R})$;

\item  $b_{0}$ is a constant advertising effectiveness factor, $b_0 \geq 0$;

\item  $b_{1}(\cdot ) $ is the density function of the
time lag between the advertising expenditure $z$ and the
corresponding effect on the goodwill level, $b_1(\cdot) \in
L^2([-R,0];\mathbb{R}_+)$;

\item  $x$ is the level of goodwill at the
beginning of the advertising campaign, $x \geq 0$;

\item  $\theta (\cdot )$ and $\delta (\cdot )$
are respectively the goodwill and the spending rate before the
beginning, $\theta(\cdot) \geq 0$, with $\theta(0)=x$, and
$\delta(\cdot) \geq 0$.

\end{itemize}

Note that when $a_{1}(\cdot )$, $b_{1}(\cdot )$ are identically
zero, equation (\ref{eq:SDDE}) reduces to the classical model
contained in the paper by Nerlove and Arrow (1962). We assume that
the goodwill and the investment in advertising at each time $s$
cannot be negative:
\begin{equation}\label{eq:vincoliadv}
\gamma(s)\ge 0, \qquad z(s) \ge 0,\quad  \forall s \in [t,T].
\end{equation}

Finally, we define the objective functional as
\begin{equation}
J(t,x;z(\cdot ))=\varphi
_{0}(\gamma(\taunew))-\int_{t}^{\taunew}h_{0}(z(s))\,ds,
\label{eq:obj-orig}
\end{equation}
where $\varphi _{0}$ is a concave utility function, $h_{0}$ is a
convex cost function, and the dynamic of $\gamma$ is determined by
(\ref{eq:SDDE}). The
functional $J$ has to be maximized over some set of admissible controls $%
\mathcal{U}$, for instance $L^{2}([t,T];\mathbb{R}_{+})$, the space
of square integrable nonnegative functions.

\section{The dynamic programming approach}
\label{sezioneDP} The Dynamic Programming (DP) approach to optimal
control problems can be summarized in four main steps (see for
instance Fleming and Rishel \cite{FlemingRishel} for the DP  in
the finite dimensional case and Li and Yong \cite{LiYong} for the
DP in the infinite dimensional case):
\begin{itemize}
\item[(i)] letting the initial data vary, calling {\it value
function} the supremum of the objective functional and writing an
equation whose candidate solution is the value function: the
so-called DP Principle, together with its infinitesimal version,
the Hamilton-Jacobi-Bellman (HJB) equation; \item[(ii)] solving
(whenever possible) the HJB equation to find the value function;
\item[(iii)] proving that the present value of the optimal control
strategy can be expressed as a function of the present value of
the optimal state trajectory: a so-called closed loop (or
feedback) relation for the optimal control; \item[(iv)] solving,
if possible, the Closed Loop Equation (CLE), i.e. the state
equation where the control is replaced by the closed loop
relation: the solution is the optimal state trajectory and the
optimal control strategy is consequently derived from the closed
loop relation.
\end{itemize}

Such method, when applicable, allows one to give a powerful
description of the optimal paths of an optimal control problem.

First of all we clarify that the two models above are not easy to
manage with the DP approach as they presents two special
difficulties.
\begin{itemize}
    \item  The state equation is a Delay Differential
Equation while the DP approach is generally formulated for
controlled Ordinary Differential Equation (ODE). One way to
approach the issue (for a different one, see e.g. Kolmanowskii and
Shaikhet \cite{kolmanowski}) is to rewrite
    the DDE as an ODE in an infinite dimensional space, which plays the role
    of the state space. We use in the sequel the techniques
    developed by Delfour, Vinter and Kwong (see Section
    \ref{sezioneriscrittura} below for explanation and
    Subsection \ref{subsliterature} for references).
    It must be noted that the resulting infinite dimensional control
    problem is harder than the ones usually treated in the literature (see e.g. \cite{LiYong})
    due to the unboundedness of the control operator and the non-analyticity of the semigroup
    involved (see again Subsection \ref{sezioneriscrittura}).

    \item Both problem feature pointwise constraints on the state
    variable, see (\ref{eq:vincoliak}), (\ref{eq:vincoliadv}). Their
    presence makes the problem much more difficult, and only a few results in special cases
    (different from the one treated here) are available in the literature. Indeed for such
    problems in infinite dimension there is no well established theory.
    This fact is at the basis of the theoretical problem contained in the paper \cite{Boucekkine}
    and mentioned in \cite{FGAK} point (II) in the introduction: show that
    the candidate optimal trajectory satisfies the pointwise
    constraints (\ref{eq:vincoliak}).
\end{itemize}
To overcome such difficulties in \cite{FGAK} we show that for our
special problem we can exhibit an explicit solution of HJB
equation. This is the key result that allows to complete the DP
approach in \cite{FGAK}.

Here, since we do not want to write the utility functions in a
fixed explicit form (like the CRRA used in
\cite{Boucekkine,FGAK}), we cannot obtain an an explicit solution
of HJB equation. Therefore we would like (here and in the future)
to perform the following steps: proving existence (and possibly,
uniqueness) for the HJB equation, then some
 theoretical results of type (iii) and (iv) above, and hopefully some subsequent
numerical approximation. This is a wide and difficult program. In
this paper we take just a first step towards the scope: existence
results for the HJB equation.


\subsection{The literature on Delay Differential Equations and on
Dynamic Programming in infinite dimensions} \label{subsliterature}
For Delay Differential Equations a recent, interesting and
accurate reference is the book by Diekmann, van Gils, Verduyn,
Lunel and Walther \cite{Diekmann}.

The idea of writing delay system using a Hilbert space setting was
first due to Delfour and Mitter \cite{DelfourMitter1},
\cite{DelfourMitter2}. Variants and improvements were proposed by
Delfour \cite{Delfour1}, \cite{Delfour2}, \cite{Delfour3}, Vinter
and Kwong \cite{VinterKwong}, Delfour and Manitius
\cite{DelfourManitius}, Ichikawa \cite{Ichikawa2} (see also the
precise systematization of the argument in chapter 4 of
Bensoussan, Da Prato, Delfour and Mitter \cite{BDDM}).

The optimal control problem in the (linear) quadratic case is
studied in Vinter, Kwong \cite{VinterKwong}, Ichikawa
\cite{Ichikawa}, Delfour, McCalla and Mitter
\cite{DelfourMcCallaMitter}. In that case the
Hamilton-Jacobi-Bellman reduces to the Riccati equation.

The study of Hamilton-Jacobi-Bellman equation in Hilbert spaces,
started with the papers of Barbu and Da Prato
\cite{BarbuDaPrato1}, \cite{BarbuDaPrato3}, \cite{BarbuDaPrato2},
is a large and diversified research field. We recall that the best
one may achieve is a ``classical'' solution of HJB equations (i.e.
solutions that are differentiable in time and state) since this
allows to get a more handleable closed loop form of the optimal
strategy. Since classical solutions are not always available,
there is a second stream in the literature that studies the
existence of ``weak'' solutions (i.e. solutions that are not
differentiable)\footnote{The most general concept of weak solution
is the one of viscosity solution, introduced by Crandall and Lions
in the finite dimensional case and then applied to infinite
dimension by the same authors, see \cite{CILuserguide} for an
introduction to the topic and further references.}. In this paper
we investigate existence of a weak-type solution (that we call
\emph{ultraweak}, see Section \ref{sezioneultraweak}) that are
limits of classical solutions.
Up to now, to our knowledge, the existence of such solutions for
the HJB equation in cases where the state equation is a Delay
Differential Equation has not been studied in the literature
(apart from the linear quadratic case). In the economic literature
the study of infinite dimensional optimal control problems that
deals with vintage/heterogeneous capital or advertising models is
a quite recent tool but of growing interest: see for instance
\cite{baruccigozzi}, \cite{FeichtingerJET}, \cite{faggian3},
\cite{GMTrento,GMS}.

\section{The state equation in an infinite dimensional setting.}
\label{sezioneriscrittura} In this section we show how to rewrite
the state equations of our examples as controlled ODE's in a
suitable Hilbert space. We do it thoroughly for the first example,
as the second is similar and simpler.

\subsection{Notation and preliminary results}
In this section we recall some general results on delay
differential equations (DDE) and on the related Hilbert space
approach, as applied to our case. The reader is referred to  the
book  by Bensoussan, Da Prato, Delfour and Mitter \cite{BDDM} for
details.  We consider from now on fixed $\ritardo>0$, and $a>0$.
With notation similar to that of \cite{BDDM}, given $\taunew>t\geq
0$ and $z \in L^2([t-\ritardo, \taunew],\mathbb{R})$ (or $z\in
L^2_{loc} ([t-\ritardo,+\infty),\mathbb{R}) $), for every $s\in
[t,\taunew]$ (or $s\in [t,+\infty)$) we call $z_s \in
L^2([-\ritardo,0];\mathbb{R})$ the function
\[
\left\lbrace
\begin{array}{ll}
z_s \colon [-\ritardo,0] \to \mathbb{R}\\
z_s(\sigma) \nd z(s+\sigma)
\end{array}
\right.
\]

Given a control $c\in\mathcal{A}$ we consider the the
following delay differential equation: \be \label{eqoriginale}
\left\{ \begin{array}{ll}
\dot{k}(s)=ak(s)-ak(s-\ritardo) - c(s) + c(s-\ritardo) \;\;\;\; for \; s\in[t,\taunew]\\
({k}(t), k_t, c_t) =(\phi^0, \phi^1, \omega)\in \mathbb{R} \times
L^2([-\ritardo,0];\mathbb{R}) \times L^2([-\ritardo,0];\mathbb{R})
\end{array} \right.\ee
where $k_t$ and $c_t$ are interpreted by means of the definition
above. Note that in the delay setting the initial data are a
triple, whose first component is the state, the second and third
are respectively the history of the state and the history of the
control up to time $t$ (more precisely, on the interval
$[t-\ritardo,t]$). The equation does not make sense pointwise,
but has to be regarded in integral sense. We give now a more
precise existence result and an estimate on the solution:
\begin{Theorem}
\label{thesistenza} Given an initial condition  $(\phi^0, \phi^1,
\omega) \in \mathbb{R} \times L^2([-\ritardo,0];\mathbb{R}) \times
L^2([-\ritardo,0];\mathbb{R})$ and a control $c \in L^2([t,
\taunew],\mathbb{R})$  there exists a unique solution $k(\cdot)$ of
(\ref{eqoriginale}) in $W^{1,2}([t, \taunew],\mathbb{R})$. Moreover
there exists a positive constant $C(\taunew-t)$ such that \be
\label{stimesoluzioneoriginaria} |k|_{W^{1,2}([t,
\taunew],\mathbb{R})} \leq C(\taunew-t) \Big ( |\phi^0|+
|\phi^1|_{L^2([-\ritardo,0];\mathbb{R})} +
|\omega|_{L^2([-\ritardo,0];\mathbb{R})} + | c|_{L^{2}([t,
\taunew],\mathbb{R})} \Big ) \ee
\end{Theorem}
\begin{proof}
See \cite{BDDM} Theorem 3.3, p.217 for the first part and Theorem
3.3 p.217, Theorem 4.1 p.222 and p.255 for the second statement.
\end{proof}
In view of the continuous embedding $W^{1,2}([t,
\taunew],\mathbb{R}) \hookrightarrow C^0([t, \taunew],\mathbb{R})$
we have also:
\begin{Corollary}
\label{stimanormainfinito} There exists a positive constant
(possibly different from the one above) $C(\taunew-t)$ such that \be
|k|_{C^0([t, \taunew],\mathbb{R})} \leq C(\taunew-t) \Big (
|\phi^0|+ |\phi^1|_{L^2([-\ritardo,0];\mathbb{R})} +
|\omega|_{L^2([-\ritardo,0];\mathbb{R})} + |c|_{L^{2}([t,
\taunew],\mathbb{R})} \Big ) \ee
\end{Corollary}

\bigskip\noindent We consider now the continuous linear application $L$ with
norm $\| L \|$
\[
\begin{array}{l}
L \colon C([-\ritardo,0],\mathbb{R}) \to \mathbb{R}\\
L \colon \varphi \mapsto \varphi(0) - \varphi(-\ritardo)
\end{array}
\]
and then define $\mathcal{L}^t$ as follows \be
\label{eqdefLstorto}
\begin{array}{ll}
\mathcal{L}^t \colon C_c([t-\ritardo,\taunew],\mathbb{R}) \to L^2([t, \taunew],\mathbb{R})\\
where \; \; \mathcal{L}^t (\phi) \colon s \mapsto L(\phi_s) \;\;\; for \; s\in[t, \taunew]\\
\end{array}
\ee where $C_c(t-\ritardo,\taunew;\mathbb{R})$ is the set of real
continuous functions having compact support contained in
$(t-\ritardo,\taunew)$
\begin{Theorem}
\label{thestensioneL} The linear operator $\mathcal{L}^t \colon
C_c([t-\ritardo,\taunew],\mathbb{R}) \to L^2([t,
\taunew],\mathbb{R})$ has a continuous extension $\mathcal{L}^t
\colon L^2([t-\ritardo,\taunew], \mathbb{R}) \to L^2([t,
\taunew],\mathbb{R})$ with norm $\leq \|L \|$ .
\end{Theorem}
\begin{proof}
See \cite{BDDM} Theorem 3.3, p. 217.
\end{proof}
\noindent Using the ``${L}$'' notation we can rewrite
(\ref{eqoriginale}) as
\[
\left\{ \begin{array}{ll}
\dot{k}(s)=a{L}(k_s)- {L}(c_s) \;\;\;\; for \; s\in[t,\taunew]\\
({k}(t), k_t, c_t)= (\phi^0, \phi^1,\omega) \in \mathbb{R} \times
L^2([-\ritardo,0];\mathbb{R}) \times L^2([-\ritardo,0];\mathbb{R})
\end{array} \right.
\]
and using the ``$\mathcal{L}^t$'' notation we can rewrite
(\ref{eqoriginale}) as \be \label{eqcome513} \left\{
\begin{array}{ll}
\dot{k}(s)=a(\mathcal{L}^tk)(s)- (\mathcal{L}^tc)(s)\;\;\;\; for \; s\in[t,\taunew]\\
({k}(t), k_t, c_t)= (\phi^0, \phi^1,\omega) \in \mathbb{R} \times
L^2([-\ritardo,0];\mathbb{R}) \times L^2([-\ritardo,0];\mathbb{R})
\end{array} \right.
\ee There follows another step towards the setting in infinite
dimension that we intend to use. So far, the history of the control
and of the trajectory were kept separated. Indeed one may note that
the delay system depends jointly on those data. Such joint
dependence is exploited in the sequel to reduce the dimension of the
state space. We then need to add some more notation to make this
more explicit.

\begin{itemize}
\item Given $u\in L^2([t-\ritardo, \taunew],\mathbb{R})$ we define
the function $e_+^t u\in L^2([t-\ritardo, \taunew],\mathbb{R})$ as
follows
\[
e_+^t u\colon [t-\ritardo,\taunew] \to \mathbb{R}, \s\s e_+^t
u(s)=\left\lbrace
\begin{array}{ll}
u(s) & s\in [t,\taunew]\\
0    & s\in [t-\ritardo,t)
\end{array}
\right.
\]
\item Given $u\in L^2([-\ritardo,0];\mathbb{R})$ we define the
function $e_-^0 u\in L^2([t-\ritardo, \taunew],\mathbb{R})$ as
follows
\[
e_-^0 u\colon [t-\ritardo,\taunew] \to \mathbb{R}, \s\s e_-^0
u(s)=\left\lbrace  \begin{array}{ll}
0 & s\in [t,\taunew]\\
u(s-t)    & s\in [t-\ritardo,t)
\end{array}
\right.
\]

\item Given a function $u\in L^2([-\ritardo,0];\mathbb{R})$  and
$s\in[t,\taunew]$ we define the function $\eta(s)u \in
L^2([-\ritardo,0];\mathbb{R})$ as follows
\[
\eta(s) u\colon [-\ritardo,0] \to \mathbb{R}, \s\s (\eta(s)
u)(\theta)= \left\lbrace
\begin{array}{ll}
u(-s+t+\theta) & \theta\geq -\ritardo+s-t\\
0 & \theta<-\ritardo+s-t
\end{array}
\right.
\]
\end{itemize}

Note that $k=e^t_+k+e^0_-\phi^1$, and $c=e^t_+c+e^0_-\omega$, then
we can separate the solution $k(s)$, $s\geq t$ and the control
$c(s)$, $s\geq t$ from initial data $\phi^1$ and $\omega$: \be
\label{eqcome514} \left\{
\begin{array}{ll}
\dot{k}=a\mathcal{L}^t e^t_+k- \mathcal{L}^te^t_+c +
a\mathcal{L}^t e^0_-\phi^1-
\mathcal{L}^te^0_-\omega \\
k(t)= \phi^0 \in \mathbb{R}
\end{array} \right.
\ee

Note that system (\ref{eqcome514}) does not directly use the
initial function $\phi^1$ and $\omega$ but only the sum of their
images $a\mathcal{L}^t e^0_+ \phi^1- \mathcal{L}^te^0_-\omega$. We
need a last step before we can write the delay equation in Hilbert
space. We introduce the operator \be \label{defdiLbar}
\left\lbrace
\begin{array}{ll}
\overline{L}\colon L^2([-\ritardo,0];\mathbb{R}) \to L^2([-\ritardo,0];\mathbb{R})\\
(\overline{L} \phi^1)(\alpha) \nd L(est(\phi^1)_{-\alpha}))
\s\s\s\s\s \alpha \in (-\ritardo,0)
\end{array}
\right. \ee where $est(\phi^1)$ is the function
$\mathbb{R}\to\mathbb{R}$ that achieves value $0$ out of
$(-\ritardo,0)$ and that is equal to $\phi^1$ in $(-\ritardo,0)$
(the same for $\omega$).

Note that the operator $\overline{L}$
    is continuous (see \cite{BDDM} page 235), moreover
\[
a\mathcal{L}^t e^0_- \phi^1 (s) - \mathcal{L}^te^0_-\omega (s)  =
 (\eta(s) (a\overline{L}\phi^1 - \overline{L} \omega)) (0) \s\s\s\s for \s s\geq
 t.
\]
Therefore, if we set \be \label{eqformadixi1}
x^1\nd(a\overline{L}\phi^1 - \overline{L} \omega), \ \
x^0\nd\phi^0,\ee then we can rewrite (\ref{eqcome514}) and
consequently (\ref{eqoriginale}) as \be \label{eqcome519} \left\{
\begin{array}{ll}
\dot{k}(s)=(a\mathcal{L}^t e^t_+k)(s) - (\mathcal{L}^te^t_+c)(s) +
(\eta(s)x^1)(0) \s\s\s\s for \s s\geq
 t \\
k(t)= x^0 \in \mathbb{R}
\end{array} \right.
\ee where $\mathbb{R} \times L^2([-\ritardo,0];\mathbb{R}) \ni x
\nd (x^0,x^1)$, $c\in \mathcal{A}$. Note that (\ref{eqcome519}) is
meaningful for all $x \in \mathbb{R}\times
    L^2([-\ritardo,0];\mathbb{R})$, also when $x^1$ is not of the form
    (\ref{eqformadixi1}). So we have embedded the original system
    (\ref{eqoriginale}) in a family of systems of the form
    (\ref{eqcome519}).

\subsection{The state equation of the AK model in the Hilbert setting}

We now work on the following Hilbert space
$$M^2 \nd \mathbb{R} \times L^2([-\ritardo,0];\mathbb{R})$$
where the scalar product between two elements
$\phi=(\phi^0,\phi^1)$ and $\xi=(\xi^0,\xi^1)$ is given by
$$\left\langle \phi, \xi \right\rangle_{M^2} \nd \left\langle
\phi^1, \xi^1 \right\rangle_{L^2} + \phi^0 \xi^0.$$ Next we
consider the homogeneous system
\[
\left\lbrace
\begin{array}{ll}
\dot z (s) = (a\mathcal{L}^0z)(s)\\
(z(0),z_0)=\phi\in M^2 
\end{array}
\right.
\]
 and define the family of continuous linear
transformations on $M^2$
\[
\left\lbrace
\begin{array}{ll}
S(s) \colon M^2 \to M^2\\
\phi \mapsto S(s) \phi \nd (z(s), z_s).
\end{array}
\right.
\]
Then $\{S(s)\}_{s\geq 0}$  is a $C_0$ semigroup on $M^2$ whose
generator is
\[
\left\lbrace
\begin{array}{ll}
D(G)=\left\lbrace (\phi^0,\phi^1)\in M^2 \s
: \s \phi^1\in W^{1,2}(-\ritardo,0) \s and\s \phi^0=\phi^1(0) \right\rbrace \\
G(\phi^0,\phi^1)=(aL\phi^1, D\phi^1)
\end{array}
\right.
\]
where $D \phi^1$ is the first derivative of $\phi^1$. A proof of
this assertion can be found in \cite{BDDM}, Chapter 4.

Note that the second component $\phi^1$ of the elements of $D(G)$
is in $C([-\ritardo,0],\mathbb{R})$ so, with a slight abuse of
notation,
 we can re-define $L$ on $D(G)$ in the following way
\[
\left\lbrace
\begin{array}{ll}
L\colon D(G) \to \mathbb{R} \\
L(\phi^0,\phi^1)=L\phi^1
\end{array}
\right.
\]
Moreover, if $D(G)$ is endowed with the graph norm, we denote with
$j$ the continuous inclusion $D(G) \hookrightarrow M^2$. Hence the
operators $G$, and $j$ are continuous from $D(G)$ into $M^2$ and
$L$ is continuous from $D(G)$ into $\mathbb{R}$. We call $G^*$,
$j^*$ and $L^*$ their adjoints, and identify $M^2$ and
$\mathbb{R}$ with their dual spaces, so that
\[
\begin{array}{ll}
G^* \colon M^2 \to D(G)'\\
j^* \colon M^2 \to D(G)'\\
L^* \colon \mathbb{R} \to D(G)'
\end{array}
\]
are linear continuous.
\begin{Definition}
The \emph{structural state} $x(s)$ at time $t\geq 0$ is defined by
\be \label{eqstatostrutturale} y(s) \nd (y^0(s),y^1(s)) \nd (k(s),
a\overline{L}(e^t_+ k)_s - \overline{L}(e^t_+ c)_s + \eta(s)x^1)
\ee
\end{Definition}\medskip\noindent
In the sequel we use $y^0$ and $y^1$ to indicate respectively the
first and the second component of the structural state. We can
give also a different, more explicit, definition: if we call
$\stackrel{\leftarrow}{k}_s, \stackrel{\leftarrow}{c}_s \in
L^2([-\ritardo,0];\mathbb{R})$ the applications
\[
\begin{array}{l}
\stackrel{\leftarrow}{k}_s \colon \theta \mapsto -k(s-\ritardo-\theta)\\
\stackrel{\leftarrow}{c}_s \colon \theta \mapsto
-c(s-\ritardo-\theta)
\end{array}
\]
the structural state can be written as \be \label{staterestated}
y(s) \nd (k(s), a \stackrel{\leftarrow}{k}_s -
\stackrel{\leftarrow}{c}_s + \eta(s)x^1). \ee

Eventually, we write  the delay equation in the Hilbert space
$M^2$ by means of the following theorem.
\begin{Theorem}
\label{thBDDM} Let $y^0(s)$ be the solution of system
(\ref{eqcome519}) for $x\in M^2$, $c\in \mathcal{A}$ and let $y(t)$
be the structural defined in (\ref{eqstatostrutturale}). Then for
each $\taunew>0$, the state $y$ is the unique solution in
\[
\left\lbrace f\in C([t,\taunew], M^2) \s : \s \frac{\ud}{\ud s}
j^* f \in L^2([t,\taunew], D(G)')  \right\rbrace
\]
to the following equation
\be
\label{eqBDDM}
\left\lbrace \begin{aligned}
        &\frac{\ud}{\ud s} y(s) = G^* y(s) + L^* c(s)\\
        &y(t)= x.
         \end{aligned} \right.
\ee
\end{Theorem}
\begin{proof}
See \cite{BDDM} Theorem 5.1 Chapter 4.
\end{proof}


\subsection{The state equation of the advertising model in the Hilbert setting}
Similar arguments can be used for the advertising model. We write
here only the results. We call $N$, $B$ the continuous linear
functionals given by
\[
\begin{array}{l}
N \colon C([-\ritardo,0]) \to \mathbb{R}\\
N \colon \varphi \mapsto a_0 \varphi(0) + \int_{-r}^{0}
\varphi(\xi ) da_{1}(\xi )
\end{array}
\]
\[
\begin{array}{l}
B \colon C([-\ritardo,0]) \to \mathbb{R}\\
B \colon \varphi \mapsto b_0 \varphi(0) + \int_{-r}^{0}
\varphi(\xi ) db_{1}(\xi )
\end{array}
\]
Let $G$ be the generator of $C_0$-semigroup defined as:
\[
\left\lbrace
\begin{array}{ll}
D(G)=\left\lbrace (\phi^0,\phi^1)\in M^2 \s : \s \phi^1\in W^{1,2}(-\ritardo,0)
 \s and\s \phi^0=\phi^1(0) \right\rbrace \\
G(\phi^0,\phi^1)=(N\phi^1, D\phi^1)
\end{array}
\right.
\]
We define $\overline{N}$ and $\overline{B}$ in the same way we defined $\overline{L}$
in equation (\ref{defdiLbar}).
So we can write the advertising model in infinite dimensional form. We obtain:
\begin{itemize}
\item The structural state in the advertising model will have the
following expression:
\[
y(t)= (y^0(s),y^1(s)) \nd (\gamma(s),
\overline{N}(e^0_+ \gamma)_s - \overline{B}(e^0_+ z)_s + \eta(s)x^1)
\]
where $x_1=\overline{N}(\theta)-\overline{B}(\delta)$.
\item The state equation becomes
\[
\left\lbrace \begin{aligned}
        &\frac{\ud}{\ud s} y(s) = G^* y(s) + B^* z(s)\\
        &y(t)= x.
         \end{aligned} \right.
\]
\end{itemize}

\section{The target functional and the HJB equation} \label{sectionHJB}

%
%
%


We now rewrite the profit functional for the first example in
abstract terms, noting that a similar reformulation holds for the
target functional of the second example. We consider a control
system governed by the linear equation described in Theorem
\ref{thBDDM}. We assume that the set of admissible controls is
defined by
\[
{\mathcal{A}} \nd \{c(\cdot)\in L^2([t, \taunew],\mathbb{R})\; : \;
c(\cdot) \geq 0\; and \; y^0(\cdot)\geq 0 \}
\]
As usual, the trajectory $y(\cdot)$ (and then $y^0(\cdot)$)
depends on the choice of the control $c(\cdot)$, and of initial
time and state, i.e. $y(\cdot)=y(\cdot; t,x, c(\cdot))$, but we
write it explicitly only when needed.

In order to apply the results contained in \cite{vincoli} and
recalled in the Appendix, we reformulate the maximization problem
as a minimization problem. At the same time we take the
constraints into account by modifying the target functional as
follows. If $h_0$ and $\phi_0$ are the concave {\it u.s.c.}
functions appearing in (\ref{target}), then we define
\[
\begin{array}{l}
h\colon \mathbb{R} \to \overline{\mathbb{R}}\\
h(c)=   \left\lbrace \begin{array}{ll}
    -h_0(c) &if \; c\ge0\\
    +\infty &if\; c<0
    \end{array} \right .
\end{array}
\]
\[
\begin{array}{l}
\phi\colon \mathbb{R} \to \overline{\mathbb{R}}\\
\phi(r)=   \left\lbrace \begin{array}{ll}
    -\phi_0(r) &if \; r\ge0\\
    +\infty &if\; r<0
    \end{array} \right .
\end{array}
\]
Moreover we set
\[
\begin{array}{l}
g\colon \mathbb{R} \to \overline{\mathbb{R}}\\
g(r)=   \left\lbrace \begin{array}{ll}
    0 &if \; r\geq 0\\
    +\infty &if\; r<0
    \end{array} \right .
\end{array}
\]
 Both $h$, $\phi$ and $g$  are convex {\it l.s.c.}
functions on $\mathbb{R}$. Then we define the {\it target
functional} as
\[
J(t,x,c(\cdot))=\int_t^\taunew e^{-\rho s} [h(c(s))+ g(y^0(s))] \ud
s + \phi(y^0(\taunew))
\]
with $c$ varying in the {\it set of admissible controls} $L^2([t,
\taunew],\mathbb{R})$. It is easy to check that the problem of
maximizing (\ref{target}) in the class $\mathcal{A}$ is equivalent
to that of minimizing $J$ on the whole space $L^2([t,
\taunew],\mathbb{R})$. Then the original maximization problem for
the AK-model has been reformulated as the following {\it abstract
minimization problem}:

 \be\label{ocp}\inf\{J(t,x,c(\cdot))\ :\ c\in L^2([t,
\taunew],\mathbb{R}),\ \textrm{and} \ y\ \textrm{satisfies}\
(\ref{eqBDDM})\},\ee Moreover, HJB equation is naturally
associated to such minimization problem by DP, and it is given by
\[
\tag{HJB} \label{HJB} \left\lbrace
\begin{array}{l}
\partial_t v(t,x) + \lla \nabla v(t,x) , G^* x \rra - \HAM(t,\nabla v(t,x))+e^{-\rho t}g(x) = 0\\
v(\taunew,x)=\phi_0(x)
\end{array}
\right .
\]
with $\HAM$ defined as follows
\[
\left\lbrace
\begin{array}{ll}
\HAM\colon [0,T] \times D(G) \to \mathbb{R}\\
\HAM(t, p) \nd \sup_{c\geq 0} \left\{ -L(p) c - e^{-\rho t}h_0(c)
\right\}=e^{-\rho t}h^*(-e^{\rho t}L(p))
\end{array}
\right .\] where $h^*$ is the Legendre transform of the convex
function $h$. We refer to $\HAM$ as to the \textit{Hamiltonian} of
the system\footnote{Note that, following the usual definition, the
Hamiltonian should be indeed $\lla p , G^* x \rra - \HAM(t,
p)+e^{-\rho t}g(x)$. Here, for commodity of notation, we put aside
of the Hamiltonian the terms which are linear or constant in $p$.}.

The abstract framework is then set, and we are ready to perform
Dynamic Programming.

\section{The value function as ultraweak solution of HJB}
\label{sezioneultraweak} We define the value function of the optimal
control problem described in the previous sections as
$$W(t,x)\nd\inf_{c(\cdot)\in L^2([t,T];\mathbb{R})}J(t,x,c(\cdot)).$$
Our objective here is to provide a suitable concept of solution of
HJB, so that the value function $V$ is a solution, in such sense.

We recall that in \cite{vincoli}
 it is shown that, if the data satisfy certain assumptions (involving convexity, semicontinuity,
 and coercivity of $h$), then
  the value function of an optimal control
problem with state constraints of type (\ref{ocp}) is indeed the
unique \textit{weak} solution to a HJB equation of type
(\ref{HJB}), as there proved and here recalled in the Appendix,
Theorem \ref{EUweek}. Note that some coercivity for the function
$h$ is indeed lacking in our case, as the prototype of $h_0$ is
$\frac{c^{1-\sigma}}{1-\sigma}$ as mentioned before, which is
sublinear on the positive real axis. This causes the Hamiltonian
of the problem - that is related to the Legendre transform of
$h_0$ - to be possibly nonregular, so that all previous definition
of solutions do not apply. (Note indeed that, as more precisely
stated in the Appendix, a weak solution is limit of strong
solutions of approximating equations, while a strong solution is
itself limit of classical solutions of approximating equations.
All of these notions require the Hamiltonian to be differentiable
with respect to the co-state variable $p$.)

 Here we are about to define a {\it ultraweak} solution as
limit of weak solutions to (\ref{HJB}). Note that the concept of
solution is indeed generalized, although not in the same direction
as before, due to the presence of possibly nonregular
Hamiltonians.

 According to the notation in \cite{Fa2}, if $X$ and $Y$
are Banach spaces, we set
\begin{equation*}\begin{split}
&Lip(X;Y)=\{f:X\to Y ~:~\lips f\lipd:=\sup_{x,y\in X,~x\neq y}
\frac{\vert f(x)-f(y)\vert_{Y}}{\vert x-y\vert_X} <+\infty\}\\
&\clip(X):=\{f\in C^1(X)~:~ \lips f^\prime\lipd<+\infty\}\\
&C_p(X,Y):=\{f:X\to\mathbb{R}~:~\vert f\vert_{C_p}:=\sup_{x\in X}
{\vert f(x)\vert_Y\over 1+\vert x\vert_X^p}<+\infty\},\ \ \
C_p(X):=C_p(X,\mathbb{R}).\\
\end{split}\end{equation*}
Moreover we set
\begin{equation*}\Sigma_0(X):=\{w\in C_2(X)\ :\ w\ {\rm is\ convex,\ } w\in\clip(X) \}\end{equation*}
\begin{equation*}\begin{split}\mathcal{Y}([0,\t]\times X)
=\{w:[0,\t]\times X\to \mathbb{R}\ :\ w\in C([0,\t],&C_2(X)),\
\\w(t,\cdot)\in&\Sigma_0(X),\ \nabla w\in C([0,\t], C_1(X,{X^\prime}))\}.\\
\end{split}
\end{equation*}

\begin{Definition}\label{UWS}
 We say that a function $V$ is a \textrm{ultraweak} solution to
\[\left\lbrace
\begin{array}{l}
\partial_t v(t,x) + \lla \nabla v(t,x) , G^* x \rra - \HAM(t,\nabla v(t,x))+e^{-\rho t}g(x) = 0\\
v(\taunew,x)=\phi_0(x)
\end{array}
\right .
\]
if there exists a sequence $\{\HAM_n\}_n$ of functions in the
space $\mathcal{Y}([0,T]\times D(G))$, such that
$\HAM_n\uparrow\HAM$ pointwise, and
$$V(t,x)=\lim_{n\to+\infty}V_n(t,x)=\inf_{n\ge0}V_n(t,x)$$
with $V_n$ the unique weak solutions to
\[\left\lbrace
\begin{array}{l}
\partial_t v(t,x) + \lla \nabla v(t,x) , G^* x \rra - \HAM_n(t,\nabla v(t,x))+e^{-\rho t}g(x) = 0\\
v(\taunew,x)=\phi_0(x)
\end{array}
\right .
\]
\end{Definition}

 Note that any weak solution $V$ is convex in the state variable
 $x$, but not necessarily ${\it l.s.c}$ in $(t,x)$. We are able to
 prove an existence result for equation $(\ref{HJB})$ by proving that the
 value function of the control problem set in the previous section
is an ultraweak solution.

\begin{Theorem} \label{EUWS} The value function $W$ of the optimal control problem (\ref{ocp})  is
an ultraweak solution of (\ref{HJB}). \end{Theorem}

\emph{Proof.} First of all we need to construct a sequence of
Hamiltonians $\HAM_n$ having the properties required by the
definition above. We choose
$$\HAM_n(t,p):=e^{-\rho t}h_n^*(-e^{\rho t}L(p))$$
with
$$h_n(c)=h(c)+\frac{1}{2n}\vert c\vert^2,\ n\in\mathbb{N}.$$
 Indeed if we denote
with $S_nf(x)=\inf_{y\in \mathbb{R}}\left\{f(y)+\frac{n}{2}\vert
x-y\vert^2\right\}$ the Yosida approximation of a function $f$,
then it is easy to check that
$[S_nf]^*(x)=f^*(x)+\frac{1}{2n}\vert x\vert^2,$ so that
$$ h_n^*(c)=S_n(h^*)(c).$$
Being $h_n^*$ the Yosida approximations of a \textit{l.s.c.}
convex function, they result to be Frech\'et differentiable with
Lipschitz gradient, with Lipschitz constant $[(h_n^*)^\prime]_L\le
n$. Moreover, as $h_n$ is a decreasing sequence, $\HAM_n$ is then
increasing, as required by Definiton \ref{UWS}. Hence the
assumptions in Theorem \ref{EUweek} are satisfied for the problem
of minimizing the functional
$$J_n(t,x,c)=J(t,x,c)+\frac{1}{2n}\int_t^\taunew e^{-\rho
s}\vert c(s)\vert^2ds$$
 in $L^2([t, \taunew],\mathbb{R})$, and we easily derive as a
 consequence the following result.
 \begin{Lemma} Let
 $$W_n(t,x)\nd \inf_{c\in L^2([t, \taunew],\mathbb{R})} J_n(t,x,c), $$
 be the value functions of the approximating optimal control
 problem. Then $W_n$ is convex in $x$ and l.s.c. in $x$ and $t$, and it is
the unique weak solution of
\[
\left\lbrace
\begin{array}{l}
\partial_t v(t,x) +\lla \nabla v(t,x) , G^* x \rra - \HAM_n( t,\nabla v(t,x))+
e^{-\rho t}g(x)= 0\\
v(\taunew,x)=\phi(x)
\end{array}
\right .
\]
Moreover there exists $c_n^*\in L^2([t, \taunew],\mathbb{R})$
optimal for the approximating problems, i.e.
$W_n(t,x)=J_n(t,x,c_n^*)$.
\end{Lemma}

To complete the proof we need to show that $W_n(t,x)\downarrow
W(t,x)$.

\begin{Lemma} The value function of (\ref{ocp}) is given by
$$W(t,x)=\lim_{n\to \infty}W_n(t,x)=\inf_nW_n(t,x).$$\end{Lemma}
{\it Proof.} By definition of $J_n$, for all
 $t$, $x$ and $n$ we have $J_n(t,x,c)\ge J_{n+1}(t,x,c)$ for all
 admissible controls $c$, so that
$$W_n(t,x)\ge W_{n+1}(t,x),$$ and $\{W_n(t,x)\}_n$ is a decreasing sequence.
As a consequence, an ultraweak solution $V$ of HJB exists, and it
is given by
$$V(t,x)\nd
\lim_{n\to\infty}W_n(t,x)=\inf_{n\in\mathbb{N}}W_n(t,x).$$ Next we
show that a solution $V$ built this way necessarily coincides with
$W$. Note that
$$J(t,x,c)\le J_n(t,x,c),\ \ \forall c\in L^2([t, \taunew],\mathbb{R}),$$ so that by
taking the infimum and then passing to limits, we obtain
\begin{equation}\label{mezzastima}W(t,x)\le V(t,x).\end{equation}
 We then prove the reverse inequality. Let $\varepsilon>0$ be arbitrarily fixed,
 and
 $c_\varepsilon$ be
an $\varepsilon$-optimal control for the problem, that is
$W(t,x)+\varepsilon>J(t,x,c_\varepsilon).$ Note that, by passing
to limits as $n\to+\infty$ in
$$V(t,x)\le W_n(t,x)\le J_n(t,x,c_\varepsilon)$$
one obtains
$$V(t,x)\le J(t,x,c_\varepsilon)<W(t,x)+\varepsilon,$$
which implies, together with (\ref{mezzastima}), the thesis.

Doing so we proved the lemma and Theorem \ref{EUWS}.

\begin{Remark}Note that we do not derive any uniqueness
result for ultraweak solutions. If for instance one tries to get
uniqueness by showing that any ultraweak solution of HJB is the
value function of a certain control problem, some difficulties
arise, due to the fact that, although $h_n^*\uparrow\mathcal{H}$
if and only if there exists some $h$ such that $h_n\downarrow h$,
in general $\mathcal{H}^*\not=h$ unless some minimax condition is
satisfied, such as
$$h= \inf_n\sup_r\{cr-h_n^*(r)\}=
\sup_r\inf_n\{cr-h_n^*(r)\}=\mathcal{H}^*,$$ which is false in
general.
\end{Remark}


\section{Appendix}
\label{sezioneappendice} In this section we recall the abstract
framework and the main results contained in \cite{Fa2} and
\cite{vincoli}, regarding strong and weak solutions of HJB.

In \cite{Fa2} and \cite{vincoli} we worked in an abstract
setting on some state space denoted with $\V$. In that setting, if $H$ is a
separable Hilbert space, $\A$ is the generator of a strongly
continuous semigroup of operators on $H$, and $V$ is the Hilbert
space
 $D(\A^*)$ endowed with the scalar product $(v|w)_V:=
(v|w)_H+(\A^*v|\A^*w)_H$, then we set $\V$ equal to its dual space
endowed with the operator norm.  The semigroup generated by $\A$
can be extended in a standard way to a semigroup
$\{e^{As}\}_{s\ge0}$ on the space $\V$, with generator $A$, a
proper extension of $\A$.
\medskip

\noindent Then we assume the state equation in $\V$ is given by
\be \label{se}
\begin{cases}y^\prime(s)=Ay(s)+Bc(s), &s\in[t,\taunew]\\
y(t)=x\in\V &\\
\end{cases}\ee
with control operator $B\in L(U,\V)$ (although $B\not\in L(U,H)$),
where $U$ is the control space and $c\in L^2([t,\taunew],U)$ the
control. Such equation may be readily expressed  in mild form as
\begin{equation}\label{mild}y(s)=e^{A(s-t)}x+\int_t^s
e^{A(s-\sigma)}Bc(\sigma)d\sigma .\ \end{equation}

\begin{Remark} The role of $\V$ in the case of
the delay equation here presented is played by the space
$D(G)^\prime$, and the role of $\A$ by the operator $G^*$.\end{Remark}
\noindent Besides, we consider a target functional $J_0$, associated to the state equation,
of type
\begin{equation}\label{J in V}J(t,x,c)
=\int_t^\t\left[g\left(s,y(s)\right)+h\left(s,c(s)\right)\right]ds+
\varphi(y(\tau))\end{equation} with $h(t,\cdot)$ real, convex,
{\it l.s.c.}, coercive, and
  $g(t,\cdot)$ and $\nu$ real, convex, and $C^1(\V)$ (respectively,
  {\it l.s.c.} in $\V$)
in the $x$ variable,  as more precisely stated in the
next sections. The problem is that of minimizing $J(t,x,\cdot)$ over the set of
admissible controls $L^2([t,\t];U)$.

\begin{Remark}
Indeed, in the applications, the target functional is rather of type
\begin{equation*}J_0(t,x,c)
=\int_t^\taunew\left[\xi\left(s,y(s)\right)+\eta\left(s,c(s)\right)\right]ds+
\nu(y(\t))\end{equation*} with $\eta(t,\cdot)$ real, convex, {\it
l.s.c.}, coercive, and
  $\xi(t,\cdot)$ and $\nu$ real, convex, and $C^1(H)$ (respectively,
  {\it l.s.c.} in $H$)
in the $x$ variable, defined on $H$, but not necessarily on $\V$. Then
we need to {\rm assume} that $\xi$ and $\nu$ allow  $C^1$
(respectively, {\it l.s.c.}) extensions $g(t,\cdot)$ and $\phi$
on the space $\V$.  The
existence of such extensions is of course a strong assumption, see
\cite{Fa2} for details and comments upon this matter.
\end{Remark}
Moreover, the value function is defined as
\begin{equation}\label{value}W(t,x)=\inf_{c\in L^2([t,\t];U)} J(t,x,c),\end{equation}
 Finally, we considered the following (backward) HJB equation
 associated to the problem set in $[0,\t]\times\V$
\begin{equation}\label{HJBb}\begin{cases}v_{t}(t,x)-\H(t,B^*\nabla v(t,x))+\langle
  A x\vert  \nabla v(t,x)\rangle+g(t,x)=0, &\\
v(\t,x)=\varphi(x),&\\
\end{cases}\end{equation}
 for all $t$ in $[0,\t]$ and $x$ in $D(A)$ (indeed for all $x$ in $\V$),
 where
$$\H(t,c)=[h(t,\cdot)]^*(-c).$$
 Note that $\H$ is  well defined only for $p$ in $V$, that is a proper subspace
of $H$, to which
 $\nabla v(t,x)$ (the spatial gradient of $v$) belongs.

With such a problem in mind, we then investigate
existence and uniqueness for the following forward HJB equation
\begin{equation}\label{HJBf}\begin{cases}\phi_t(t,x)+F(t,\nabla\phi(t,x))-\langle
  Ax,\nabla\phi(t,x)\rangle=g(\t-t,x), &(t,x)\in[0,\t]\times {V^\prime}\\
\phi(0,x)=\varphi(x).&\\
\end{cases}\end{equation}
Note in fact that such a HJB is the forward version of (\ref{HJBb}) if we set
 \begin{equation*}F(t,p):=\H(t,B^*p)=
\sup_{c\in U}\{\big( -Bc\vert p\big)_U-h(t,c)\}.\end{equation*}

\subsection{Regular data and strong solutions of HJB equations.}
We first treat the case of regular data, from which the notion of
strong solution originates.

\begin{ass}\label{asst}

\begin{enumerate}

\item[1.] $A:D(A)\subset\V\to\V$ is the infinitesimal generator of
a strongly continuous semigroup $\{e^{s A}\}_{s\ge0}$ on $\V$;

\item[2.] $B\in L(U,\V)$;

\item[3.] there exists $\omega>0$ such that $\vert e^{\tau
A}x\vertv\le M e^{\omega \tau}\vert x\vertv,~\forall \tau\ge0$;

\item[4.] $F\in \mathcal{Y}([0,\t]\times V)$, $F(t,0)=0$,
$\sup_{t\in [0,\t]}[ F_p(t,\cdot)]_L<+\infty$;

\item[5.] $g\in \mathcal{Y}([0,\t]\times \V)$, $t\mapsto\lips
g_x(t,\cdot)\lipd\in L^1(0,\t)$

\item[6.] $\varphi\in\Sigma_0(\V)$;

\item[7.] $h(t,\cdot)$ is  convex, lower semi--continuous,
$\partial_c h(t,\cdot)$  is  injective for all $t\in[0,\t]$.

\item[8.] $\H\in \mathcal{Y}([0,\t]\times U)$, $\H(t,0)=0$, and
$\sup_{t\in[0,\t]}\lips \H_c(t,\cdot)\lipd<+\infty.$

\end{enumerate}\end{ass}

\begin{Definition}  Let  Assumptions \ref{asst} be satisfied. We say
that   $\phi\in C([0,\t],C_2(\V))$ is a {\rm strong} solution of
$(\ref{HJBf})$ if there exists a family $\{\phi^\e\}_\e\subset
C([0,\t],C_2(\V))$ such that:

$(i)$ $\phi^\e(t,\cdot)\in \clip(\V)$ and $\phi^\e(t,\cdot)$ is
convex for all $t\in[0,\t]$; $\phi^\e(0,x)=\varphi(x)$ for all
$x\in\V$.

$(ii)$ there exist  constants $\Gamma_1,\Gamma_2>0$ such that
$$\sup_{t\in[0,\t]}\lips\nabla\phi^\e(t)\lipd\le \Gamma_1,~
\sup_{t\in[0,\t]}\vert\nabla\phi^\e(t,0)\vert_V\le
\Gamma_2,~\forall \e>0;$$

$(iii)$  for all $x\in D(A)$, $t\mapsto\phi^\e(t,x)$ is
continuously
 differentiable;

$(iv)$ $\phi^\e\to\phi$, as $\e\to 0+$,
 in $C([0,\t],C_2(\V))$;

$(v)$ there exists  $g_\e\in C([0,\t];C_2(\V))$ such that, for all
$t\in[0,\t]$ and $x\in D(A)$,
$$\phi_t^\e(t,x)-F(t,\nabla\phi^\e(t,x))+\langle Ax,\nabla\phi^\e(t,x)\ranglev=g_\e(\t-t,x)$$
with $g_\e(t,x)\to g_0(t,x)$, and $\int_0^\t\vert
g_\e(s)-g_0(s)\vert_{C_2}ds\to 0$, as $\e\to 0+.$ \end{Definition}

The main result contained in \cite{Fa2} is the following.

\begin{Theorem}\label{exun}
Let Assumptions \ref{asst} be satisfied. There exists a unique
strong solution $\phi$ of $(\ref{HJBf})$ in the class
$C([0,T],C_2(\V))$ with the following properties:

$(i)$ for all $x\in D(A)$, $\phi(\cdot,x)$ is Lipschitz
continuous;

$(ii)$ $\phi(t,\cdot)\in\Sigma_0({V^\prime})$, for all
$t\in[0,T]$.
\end{Theorem}

Regarding applications to the optimal control problem, in
\cite{Fa3} we were able to prove what follows.

\begin{Theorem}\label{verif}
Let Assumptions \ref{asst} be satisfied, with
$F(t,p):=\H(t,B^*p)$. Let $W$ be the value function of the
control problem, and let $\phi$ be the strong
  solution of $(\ref{HJBf})$ described in Theorem \ref{exun}. Then
  $$W(t,x)=\phi(\t-t,x),~\forall t\in[0,\t],~\forall x\in\V,$$
  that is,
  the value function $W$ of the optimal control problem is the
  unique strong solution of the backward HJB equation $(\ref{HJBb})$.
\end{Theorem}
\bigskip

\subsection{Semicontinuous data and weak solutions of HJB equations.}
We then treat the case of merely semicontinuous data, from which the notion of
{\it weak} solution originates.

\begin{ass}\label{asst2} If $K$ is a convex closed subset of $\V$, we define
\begin{equation*}\Sigma_K\equiv\Sigma_K(\V):=\{\phi:\V\to(-\infty,+\infty]\
:\ \phi\ {\rm is\ convex\ and}\ l.s.c.,\ K\subset
D(\phi)\}\end{equation*} where $D(\phi)=\{x\in\V\ :\
\phi(x)<+\infty\}$, and assume:
\begin{enumerate}

\item[1.] $C:D(C)\subset\V\to\V$ is the infinitesimal generator of a strongly
continuous semigroup $\{e^{s A}\}_{s\ge0}$ on $\V$;

\item[2.] $B\in L(U,\V)$;

\item[3.] there exists
$\omega>0$ such that $\vert e^{s C}x\vertv\le e^{\omega s}\vert
x\vertv,~\forall s\ge0$;

\item[4.] $F\in \mathcal{Y}([0,\t]\times V)$, $F(t,0)=0$, $\sup_{t\in [0,\t]}[
F_p(t,\cdot)]_L<+\infty$;

\item[5.]
$g(t,\cdot)\in\Sigma_K(\V)$, for all $t\in [0,\t]$; $g(\cdot,x)$
{\it l.s.c.} and $L^1(0,\t)$ for all $x\in\V$;

\item[6.] $\varphi\in\Sigma_K(\V)$;

\item[7.]
$h(t,\cdot)$ is  convex, lower semi--continuous, $\partial_c
h(t,\cdot)$  is  injective for all $t\in[0,\t]$; moreover
$h(t,c)\ge a(t)\vert c\vert^2_U+b(t)$,  with $a(t)\ge A(T)>0$,
$b\in L^1(0,T;\mathbb{R})$;

\item[8.]
$\H\in \mathcal{Y}([0,\t]\times U)$, $\H(t,0)=0$, and
$\sup_{t\in[0,\t]}\lips \H_c(t,\cdot)\lipd<+\infty.$

\end{enumerate}
\end{ass}

\begin{Definition}\label{weak sol} Let $K\subset \V$ be a closed
  convex set, and
let $\varphi\in\Sigma_K$ and
$g(t,\cdot)\in \Sigma_K$ for all $t$ in $[0,\t]$.
Then $\phi:[0,\t]\times\V\to(-\infty,+\infty]$ is a {\rm weak}
solution of $(HJB)$ if:

 $(i)$
$\phi(t,\cdot)\in\Sigma_K,\ \forall t\in[0,\t]$;

$(ii)$ there exist sequences $\{\varphi_n\}_n\subset\Sigma_0$, and
$\{g_n\}\subset\mathcal{Y}([0,T]\times \V)$, such that
 $$\varphi_n(x)\uparrow\varphi(x), \ g_n(t,x)\uparrow
 g(t,x),\
 \forall x\in \V, \ \forall t\in[0,\t], \ as\ n\to+\infty,$$
and moreover, if $\phi_n$ is the unique strong solution of
$$\begin{cases}\phi_t(t,x)+ F(t,\nabla\phi(t,x))-\langle
  Ax,\nabla\phi(t,x)\ranglev= g_n(t,x) &(t,x)\in[0,\t]\times {V^\prime}\cr
\phi(0,x)=\varphi_n(x)&\cr\end{cases}$$ in $C([0,\t],C_2(\V))$,
then
$$\phi_n(t,x)\uparrow\phi(t,x),\ \forall
(t,x)\in[0,T]\times\V.$$
  \end{Definition}

 \begin{Remark} Since strong solution were proved in \cite{Fa2}
to be Lipschtiz  with respect to the time variable
 and $C^1$ with respect to the space variable, and the weak
 solution $\phi$ is a sup--envelop of strong solutions $\phi_n$,
then $\phi$ is lower semi-continuous
 in $[0,\t]\times\V$. For the same reason
$\phi_n$ convex in the $x$ variable implies that $\phi$ is convex in $x$ as well.
\end{Remark}

\begin{Remark}Note that the role of the convex set $K$ is played in
  the first example by the set
$$K\nd cl_{\V}(\{(x_0,x_1)\ :\ x_0\ge0\})$$
\end{Remark}

\begin{Theorem}\label{EUweek} Let Assumptions \ref{asst2} be satisfied.
Let also $g$ and $h$ be of the following type
$$g(t,x)=e^{-\rho t}g_0(x), \ \ h(t,c)=e^{-\rho t}h_0(x)$$.
Then the following properties are equivalent:

$(i)$ there exists a unique weak solution of (\ref{HJBf});

$(ii)$ At each $(t,x)\in [0,\t]\times K$ there exists an admissible control.

\noindent Moreover if $(i)$ or $(ii)$ holds, there exists an optimal pair $(c^*,y^*)$ and
$$\phi(\t-t,x)=J(t,x,c^*).$$ \end{Theorem}

\bibliographystyle{plain}
\bibliography{VintageCapital-InfiniteDimen291104}

\end{document}